\documentclass[11pt]{article}
\usepackage[latin1]{inputenc}
\usepackage{epsfig}
\usepackage{color}
\usepackage[british,english]{babel}
\usepackage{amsthm}
\usepackage{amsmath}
\usepackage{amsfonts}
\usepackage{amssymb}
\usepackage{graphicx}
\newcommand\dint{\displaystyle\int}

\newcommand{\blll}     {\begin{array}{lll}}
\newcommand{\brcl}     {\begin{array}{rcl}}
\newcommand{\barr}     {\begin{array}}
\newcommand{\earr}     {\end{array}}
\newcommand{\R}{I\!\!R}

\def\eps{\varepsilon}
\def\dsp{\displaystyle}
\def\Vronde{\mathcal{V}}
\def\P{\mathcal{P}}
\def\O{\mathcal{O}}
\def\M{\mathcal{M}}
\def\VV{\mathcal{V}}
\def\RR{\mathcal{R}}

\def\K{\mathcal K}

\def\PB{\mathcal P}
\def\R{\mathbb{R}}
\newtheorem{theorem}{Theorem}
\newtheorem{corollary}{Corollary}

\title{\bf A fast reconstruction algorithm for geometric inverse problems using topological sensitivity analysis and Dirichlet-Neumann cost functional approach}
\author{Maatoug Hassine$^1$ and Imen Kallel$^2$\\
\small $^{1}$ FSM, Monastir University, Tunisie\\
\small $^{2}$ Monastir University, Tunisie} 
\date{}
\begin{document}
\maketitle
 \noindent {\bf Abstract.} This paper is concerned with the detection of objects immersed in anisotropic media from boundary measurements. We propose an accurate approach based on the Kohn-Vogelius formulation and the topological sensitivity analysis method. The inverse problem is formulated as a topology optimization one minimizing an energy like functional. A topological asymptotic expansion is derived for the anisotropic Laplace operator. The unknown object is reconstructed using a level-set curve of the topological gradient. The efficiency and accuracy of the proposed algorithm are illustrated by some numerical results.



\def\eps{\varepsilon}
\def\dsp{\displaystyle}
\def\Vronde{\mathcal{V}}
\def\P{\mathcal{P}}
\def\O{\mathcal{O}}
\def\M{\mathcal{M}}
\def\VV{\mathcal{V}}
\def\RR{\mathcal{R}}
\def\PB{\mathcal P}
\def\R{\mathbb{R}}
\def\K{\mathcal K}



\section{Introduction}


In this work we are interested on the detection of objects immersed in an anisotropic media from overdetermined boundary data. More precisely, let $\Omega\subset \R^d,\,\, d=2,3$ be a  bounded domain with smooth boundary $\Gamma$. We assume that $\Gamma$ is partitioned into two parts $\Gamma_a$ (accessible) and $\Gamma_i$ (inaccessible) having both non-vanishing measure. 

Let $\O$ be an unknown object immersed inside the background domain $\Omega$ and having a smooth boundary $\Sigma = \partial \O$. The geometric inverse problem that we consider here can be formulated as follows:

$-$ Giving two boundary data on $\Gamma_a$; an imposed flux $\Phi\in H^{-1/2}(\Gamma_a)$ and a measured datum $\psi_m\in H^{1/2}(\Gamma_a)$.

$-$ Find the unknown location of the object $\O$ inside the domain $\Omega$ such that the solution $\psi$ of the anisotropic Laplace equation satisfies the following overdetermined boundary value problem
\begin{equation}
 \left \lbrace
\begin{array}{rll}
-\mbox{ div } (\gamma(x) \nabla \psi)&=F  & \mbox{ in } \Omega\backslash \overline{\O},\\
\gamma(x)\nabla \psi \cdot{\bf n }&=\Phi & \mbox{ on } \Gamma_a,\\
\psi &=\psi_m & \mbox{ on } \Gamma_a,\\
\psi &=0 & \mbox{ on } \Gamma_i,\\
\gamma(x) \nabla\psi \cdot{\bf n } &= 0 & \mbox{ on } \Sigma,
\end{array}
\right.
\end{equation}
where $\gamma$ is a scalar smooth function (of class $\mathcal C^1$) describing the physical properties of the medium $\Omega$, ${\bf n }$ is the exterior unit normal vector and $F\in L^2(\Omega)$ is a given source term. We assume that there exist two constants $c_0>0$ and $c_1>0$ such that $c_0\leq \gamma(x)\leq c_1,\,\forall x\in \Omega$.

In this formulation the domain $\Omega\backslash \overline{\O}$ is unknown since the free boundary $\Sigma$ is unknown. This problem is ill-posed in the sense of Hadamard. The majority of works dealing with this kind of problems fall into the category of shape optimization and based on the shape differentiation technics. It is proved in \cite{BCD} that the studied inverse problems, treated as a shape optimization problems, are severely ill-posed (i.e. unstable), for both Dirichlet and Neumann conditions on the boundary $\Sigma$. Thus they have to use some regularization methods to solve them numerically.

To solve this inverse problem, we extend the topological sensitivity analysis notion to the  anisotropic case and we suggest an alternative approach  based on the Kohn-Vogelius formulation \cite{ABB} and the topological gradient method \cite{AHM,AK,BHJM,GH,GB,HJM}. We combine here the advantages of the Kohn-Vogelius formulation as a self regularization technique and the topological gradient approach as an accurate and fast method.

In the first part of this paper, we derive a topological sensitivity analysis for a Kohn-Vogelius type functional valid for an arbitrary shaped geometric perturbation. The obtained results are based on a preliminary estimate describing the perturbation caused by the presence of a small geometry modification of the background domain $\Omega$.  The proposed mathematical analysis is general and can be adapted for various partial differential equations (PDE).

The second part of this paper is concerned with some numerical investigations. The obtained topological sensitivity function is used to built a simple, fast and accurate geometry reconstruction algorithm. The efficiency of the proposed algorithm is illustrated by some numerical results.


This work is motivated by many practical problems for which it is necessary to detect the electrical properties of a media from boundary measurements. This kind of studies was realyzed for the clinical applications such as electrical impedance tomography \cite{CIN}, the geophysical applications such as detection of the mineral deposits location in the earth \cite{parker},  industrial applications such as non-destructive testing \cite{bowler}, ... etc.

The rest of the paper is organized as follows. In Section 2, we formulate  the considered inverse problem as a topological optimization one. In Section 3, we discuss the perturbation caused by the presence of a small object inside the background domain $\Omega$. The Section 4 is devoted to the Kohn-Vogelius type function variation.  In Section 5, we derive a topological sensitivity analysis for the anisotropic Laplace operator. In Section 6, we propose a one-iteration  shape reconstruction algorithm. 
 
\section{Formulation of the inverse problem}
In this section, we give the main steps of our analysis. Firstly, we introduce the Kohn-Vogelius formulation and we define the cost function to be minimized. Secondly, we present the perturbed problems and we describe the quantity to be estimated. 

\subsection{The Kohn-Vogelius formulation}
The Kohn-Vogelius formulation rephrases the considered inverse problem into a topological optimization one. In fact, the Kohn-Vogelius formulation leads to define for any given  domain $\O\subset \Omega $ two forward problems. The first one is associated to the Neumann datum $\Phi$, it will be named the ``Neumann problem'':  
 \begin{equation}
\label{neum-pb}
 (\PB_N)\left\{
\begin{array}{l}
\mbox{Find } \psi_N \in H^1(\Omega\backslash \overline{\O}) \mbox{ solving }\\
\begin{array}{rll}
\displaystyle
- \mbox{ div } (\gamma(x) \nabla \psi_N)  =&  F &\hbox{ in }  \Omega \backslash \overline{\O}\\
\gamma(x) \nabla \psi_N \cdot{\bf n } =& \Phi &\hbox{ on } \Gamma_a \\
\psi_N =& 0 &\hbox{ on } \Gamma_i \\
\gamma(x) \nabla \psi_N \cdot{\bf n } =& 0 &\hbox{ on } \Sigma.
\end{array}
\end{array}
\right.
\end{equation}
The second one is associated to the Dirichlet datum (measured) $\psi_m$
 \begin{equation}
\label{dir-pb}
(\PB_D) \left\{
\begin{array}{l}
\mbox{Find } \psi_D \in H^1(\Omega\backslash \overline{\O}) \mbox{ solving }\\
\begin{array}{rll}
\displaystyle
- \mbox{ div } (\gamma(x) \nabla \psi_D)  =&  F &\hbox{ in }  \Omega \backslash \overline{\O}\\
\psi_D =& \psi_m &\hbox{ on } \Gamma_a\\
\psi_D =& 0 &\hbox{ on } \Gamma_i \\
\gamma(x) \nabla \psi_D \cdot{\bf n }=& 0 &\hbox{ on }  \Sigma.
\end{array}
\end{array}
\right.
\end{equation}

One can see that if $\Sigma$ coincides with the actual  boundary $\Sigma^*$ then the misfit between the solutions vanishes, ${\psi_D} = {\psi_N}$. According to this observation, we propose an identification process based on the minimization of the following energy like functional \cite{ABB}
\begin{eqnarray*}\label{misfit-func}
\K (\Omega \backslash \overline{\O}) = \int_{\Omega \backslash \overline{\O}}\gamma(x)|\nabla \psi_D - \nabla \psi_N|^2 dx.
\end{eqnarray*}
The inverse problem can be formulated as a topological optimization one
$$
\min_{\O \subset \Omega} \K (\Omega \backslash \overline{\O}).
$$
To solve this  problem we will use the topological sensitivity analysis method. 


\subsection{The sensitivity analysis method}
 This  method consists in studying the variation of  the function $\K$ with respect to the presence of a small object inside the background domain  $\Omega$.

To present the main idea of this method, we consider the case in which $\Omega$ contains a small object $\O_{z,\eps}$ that is centred at $z\in \Omega$ and has the shape $\O_{z,\eps}=z+\eps \omega\subset \Omega$, where $\eps>0$ and $\omega\subset \R^d$ is a given, fixed and bounded domain containing the origin, whose boundary $\partial \omega$ is of $\mathcal{C}^1$. The topological sensitivity analysis  leads to an asymptotic expansion of the variation $\K (\Omega \backslash \overline{\O_{z,\eps}})-\K (\Omega)$ with respect to $\eps$.

Using the Kohn-Vogelius formulation, one can define for each arbitrary location of $\O_{z,\eps}$ in the  domain $\Omega$, two forward problems. 

The first one is associated to the Neuman datum $\Phi$ and called ``the perturbed Neumann problem'' 
\begin{equation}
\label{npb}
 (\PB^\eps_N)\left\{
\begin{array}{rll}
\mbox{find } \psi^\eps_N \in H^1(\Omega \backslash \overline{\O_{z,\eps}})& \mbox{such that}\\
\displaystyle
- \mbox{ div } (\gamma(x) \nabla \psi^\eps_N)  =&  F &\hbox{ in } \Omega \backslash \overline{\O_{z,\eps}} \\
\gamma(x) \nabla \psi^\eps_N \cdot{\bf n } =& \Phi &\hbox{ on } \Gamma_a \\
\psi^\eps_N =& 0 &\hbox{ on } \Gamma_i \\
\gamma(x) \nabla \psi^\eps_N \cdot{\bf n }=& 0 &\hbox{ on }\partial\O_{z,\eps},
\end{array}
\right.
\end{equation}

The second one is associated to the Dirichlet datum $\psi_m$ and called ``the perturbed Dirichlet problem''  
 \begin{equation}
\label{dpb}
(\PB^\eps_D) \left\{
\begin{array}{rll}
\mbox{find } \psi^\eps_D \in H^1(\Omega \backslash \overline{\O_{z,\eps}})& \mbox{such that}\\
\displaystyle
- \mbox{ div } (\gamma(x) \nabla \psi^\eps_D)  =&  F &\hbox{ in } \Omega \backslash \overline{\O_{z,\eps}}\\
\psi^\eps_D =& \psi_m &\hbox{ on } \Gamma_a\\
\psi^\eps_D =& 0 &\hbox{ on } \Gamma_i \\
\gamma(x) \nabla \psi^\eps_D \cdot{\bf n }=& 0 &\hbox{ on } \partial\O_{z,\eps}.
\end{array}
\right.
\end{equation}
In order to describe the presence of the object $\O_{z,\eps}$ inside the domain $\Omega$, we will use the shape function
\begin{eqnarray*}
\K (\Omega \backslash \overline{\O_{z,\eps}}) =\int_{\Omega \backslash \overline{\O_{z,\eps}}}\gamma(x)|\nabla \psi^\eps_D - \nabla \psi^\eps_N|^2\, dx,
\end{eqnarray*}
Next, we will derive a topological sensitivity analysis for the function $\K$ with respect to the insertion of a small object $\O_{z,\eps}$ in $\Omega$. It leads to an asymptotic expansion of the form
\begin{eqnarray*}
\K (\Omega \backslash \overline{\O_{z,\eps}})=\K (\Omega)+ \rho(\eps) \delta \K(z) + o(\rho(\eps)),\quad \forall z\in \Omega,
\end{eqnarray*}
where $\eps \longmapsto\rho(\eps)$ is a scalar positive function going to zero with $\eps$. The function $z\longmapsto \delta \K(z)$ is called the topological gradient and play the role of leading term of the variation $\K (\Omega \backslash \overline{\O_{z,\eps}})-\K (\Omega)$. In order to minimize the shape function $\K$, the best location $z$ of the object $\O_{z,\eps}$ in $\Omega$  is where $\delta\K$ is most negative.

 We start our analysis by estimating the perturbation caused by the presence of the small object $\O_{z,\eps}$ in the background domain $\Omega$. We will establish in the next section two estimates describing the behavior of the perturbed solutions with respect to $\eps$. In Section \ref{cost-funct-variation}, we will discuss the cost function variation. Based on the obtained estimates, we will derive in Section \ref{asymptotic} a topological asymptotic expansion for the anisotropic laplace operator.    

\section{Estimate of the perturbed solutions}\label{sec-estimate}
In this paragraph, we establish two estimates describing the perturbation caused by the presence of the geometry modification $\O_{z,\eps}$ on the solutions of the Dirichlet and Neumann problems. To this end, we introduce two auxiliaries problems. \\
The first one is related to the Neumann problem:
 \begin{equation}
\label{ext-neum}
\left\{
\begin{array}{l}
\mbox{Find } \varphi_N \in W^1(\mathbb{R}^d\backslash \overline{\omega}) \mbox{ such that }\\
\begin{array}{rll}
\displaystyle
- \Delta \varphi_N  =&  0 &\hbox{ in }  \mathbb{R}^d \backslash \overline{\omega}\\
\varphi_N =& 0 &\hbox{ at } \infty \\
\nabla \varphi_N \cdot{\bf n } =& -  \nabla \psi_N (z)\cdot{\bf n }&\hbox{ on } \partial \omega .
\end{array}
\end{array}
\right.
\end{equation}
The second one is related to the Dirichlet problem:
 \begin{equation}
\label{ext-neum}
\left\{
\begin{array}{l}
\mbox{Find } \varphi_D \in W^1(\mathbb{R}^d\backslash \overline{\omega}) \mbox{ such that }\\
\begin{array}{rll}
\displaystyle
- \Delta \varphi_D  =&  0 &\hbox{ in }  \mathbb{R}^d \backslash \overline{\omega}\\
\varphi_D =& 0 &\hbox{ at } \infty \\
\nabla \varphi_D \cdot{\bf n } =& -  \nabla \psi_D (z) \cdot{\bf n }&\hbox{ on } \partial \omega .
\end{array}
\end{array}
\right.
\end{equation}
The functions $\varphi_N$ and $\varphi_N$ can be expressed by a single layer potential on $\partial \omega$ (see \cite{DL}) on the following way
\begin{eqnarray*}
\varphi_N(y)&=& \displaystyle \int_{\partial \omega} E(y-x)\eta_N(x)ds(x),\forall y \in \mathbb{R}^d \backslash \overline{\omega}\\
\varphi_D(y)&=& \displaystyle \int_{\partial \omega} E(y-x)\eta_D(x)ds(x),\, \forall y \in \mathbb{R}^d \backslash \overline{\omega},
\end{eqnarray*}
where $E$ is the fundamental solution of the Laplace problem in $\mathbb{R}^d$; 
$$E(y)=
\left\{
\begin{array}{ll}
-\displaystyle \frac{1}{2\pi} \log(|y|) &\mbox{ if }\,\,d=2,\\
\displaystyle \frac{1}{4\pi}\frac{1}{|y|}&\mbox{ if }\,\,d=3.
\end{array}
\right.
$$
\noindent Here $\eta_N$ and $\eta_D$ belong to $H^{-1/2}(\partial \omega)$ and solve the following integral equations  \cite{DL}
\begin{eqnarray}\label{integ-N}
&& -\frac{\eta_N(y)}{2}+\int_{\partial\omega}\nabla E(y-x) \cdot {\bf n }\,\eta_N(x)ds(x)=-\nabla \psi^0_N(z) \cdot {\bf n },\quad y\in \partial \omega \\\label{integ-D}
&&-\frac{\eta_D(y)}{2}+\int_{\partial\omega}\nabla E(y-x) \cdot {\bf n }\,\eta_D(x)ds(x)=-\nabla \psi^0_D(z) \cdot {\bf n },\quad y\in \partial \omega.
\end{eqnarray}
The Neumann  and Dirichlet  perturbed solutions  satisfy the following estimates.

\begin{theorem} \label{est-pert-sol} There exist positive constants $c>0$, independent of $\eps$, such that 
\begin{eqnarray*}
\left\| \psi^\eps_N-\psi^0_N-\eps\varphi_N((x-z)/\eps) \right\|_{1,\Omega_{z,\eps}}\leq c\eps^{d/2},\\
\\
\left\| \psi^\eps_D-\psi^0_D-\eps\varphi_D((x-z)/\eps) \right\|_{1,\Omega_{z,\eps}}\leq c\eps^{d/2}.
\end{eqnarray*}
\end{theorem}

\noindent {Proof:} In order to prove the estimates established for the Neumann and Dirichlet perturbed solutions, we consider the following generic problem
\begin{equation}\label{auxil-eps}
\left\{
\begin{array}{l}
\mbox{Find }\psi^\eps \in H^1(\Omega_{z,\eps}) \mbox{ such that }\\
\begin{array}{rll}
\displaystyle
- \mbox{ div } (\gamma(x) \nabla \psi^\eps)  =&  F &\hbox{ in } \Omega_{z,\eps} \\
\gamma(x) \nabla \psi^\eps{\cdot n } =& \Phi &\hbox{ on } \Gamma_n \\
\psi^\eps =& \psi_m &\hbox{ on } \Gamma_d \\
\psi^\eps =& 0 &\hbox{ on } \Gamma_i \\
\gamma(x) \nabla \psi^\eps {\cdot n }=& 0 &\hbox{ on }\partial\O_{z,\eps},
\end{array}
\end{array}
\right.
\end{equation}
where $\Gamma_n$ and $\Gamma_d$ are two parts of the boundary $\Gamma_a$ such that $\overline{\Gamma_a}=\overline{\Gamma_n}\cup\overline{\Gamma_d}$ and $\Gamma_n\cap\Gamma_d=\emptyset$.

\noindent One can remark here that the considered problem (\ref{auxil-eps}) has a general forme valid for the Neumann and Dirichlet cases. In fact,  if $\Gamma_n=\emptyset$, we have $\Gamma_d=\Gamma_a$ and $\psi^\eps$ solves the Dirichlet problem $ (\PB^\eps_D)$. If $\Gamma_d=\emptyset$, we have $\Gamma_n=\Gamma_a$ and $\psi^\eps$ solves the Neumann problem $ (\PB^\eps_N)$. 

In the absence of any geometry perturbation (i.e. $\eps=0$), we have $\Omega_{z,\eps} =\Omega$ and $\psi^0$ solves
\begin{equation}\label{auxil-0}
\left\{
\begin{array}{rll}
\displaystyle
- \mbox{ div } (\gamma(x) \nabla \psi^0)  =&  F &\hbox{ in } \Omega \\
\gamma(x) \nabla \psi^0 \cdot {\bf n } =& \Phi &\hbox{ on } \Gamma_n \\
\psi^0 =& \psi_m &\hbox{ on } \Gamma_d \\
\psi^0 =& 0 &\hbox{ on } \Gamma_i.
\end{array}
\right.
\end{equation}
We denote by $\varphi$ the solution to the associated exterior problem
 \begin{equation}
\label{auxil-ext}
\left\{
\begin{array}{rll}
\displaystyle
- \Delta \varphi  =&  0 &\hbox{ in }  \mathbb{R}^d \backslash \overline{\omega}\\
\varphi =& 0 &\hbox{ at } \infty \\
\nabla \varphi \cdot{\bf n } =& -  \nabla \psi^0 (z)\cdot{\bf n }&\hbox{ on } \partial \omega .
\end{array}
\right.
\end{equation}
Combining (\ref{auxil-eps}), (\ref{auxil-0}) and (\ref{auxil-ext}), one can deduce that $\phi^\eps=\psi^\eps-\psi^0-\eps \varphi((x-z)/\eps)$ is solution to 
\begin{equation}
\left\{
\begin{array}{rll}
\displaystyle
- \mbox{ div } (\gamma(x) \nabla \phi^\eps)  =&  \nabla\gamma(x).\nabla_y \varphi((x-z)/\eps) &\hbox{ in } \Omega_{z,\eps} \\
\gamma(x) \nabla \phi^\eps \cdot{\bf n } =& -\gamma(x) \nabla_y \varphi((x-z)/\eps)\cdot{\bf n } &\hbox{ on } \Gamma_n \\
\phi^\eps =& -\eps \varphi((x-z)/\eps) &\hbox{ on } \Gamma_d\cup\Gamma_i \\
\gamma(x) \nabla \phi^\eps \cdot{\bf n }=& -\gamma(x) [\nabla \psi^0- \nabla \psi^0(z)]\cdot{\bf n } &\hbox{ on }\partial\O_{z,\eps}.
\end{array}
\right.
\end{equation}
Due to the smoothness of $\gamma$ in $\Omega$, there exists $c>0$ such that
\begin{eqnarray*}
\left\| \nabla\gamma(x).\nabla_y \varphi ((x-z)/\eps) \right\|_{0,\Omega_{z,\eps}} &\leq c & \left\|\nabla_y \varphi ((x-z)/\eps) \right\|_{0,\Omega_{z,\eps}},\,\\
\left\| \gamma(x) \nabla_y \varphi((x-z)/\eps)\cdot{\bf n }\right\|_{-1/2,\Gamma_n} &\leq c & \left\| \nabla_y \varphi((x-z)/\eps) \cdot{\bf n }\right\|_{-1/2,\Gamma_n}.
\end{eqnarray*}
By trace Theorem, it follows 
\begin{eqnarray*}
\left\| \varphi((x-z)/\eps)\right\|_{1/2,\Gamma_d\cup\Gamma_i}+\left\| \nabla_y \varphi((x-z)/\eps)\cdot{\bf n }\right\|_{-1/2,\Gamma_n} \leq c\,\left\| \varphi((x-z)/\eps)\right\|_{1,\Omega_R},
\end{eqnarray*}
where $\Omega_R=\Omega\backslash\overline{B(z,R)}$, with $R>0$ is a given radius such that $\overline{\O_{z,\eps}} \subset B(z,R)$ and $\overline{B(z,R)}\subset \Omega$.

\noindent It is easy to check that the function $\varphi^\eps(x)=\varphi((x-z)/\eps)$ solves the following problem
 \begin{equation}
\label{auxil-ext-eps}
\left\{
\begin{array}{rll}
\displaystyle
- \Delta \varphi^\eps  =&  0 &\hbox{ in }  \mathbb{R}^d \backslash \overline{\O_{z,\eps}}\\
\varphi^\eps =& 0 &\hbox{ at } \infty \\
\nabla \varphi^\eps \cdot{\bf n } =& -  \nabla \psi^0 (z)\cdot{\bf n }&\hbox{ on } \partial \O_{z,\eps}.
\end{array}
\right.
\end{equation}
Since $\displaystyle \int_{\partial \O_{z,\eps}} \nabla \psi^0 (z)\cdot{\bf n }ds=0$, using the change of variable $x=z+\eps y$, one can prove that there exists $c>0$ (independent of $\eps$) such that 
\begin{eqnarray}\label{est-1}
\left\|\nabla_y \varphi^\eps \right\|_{0,\Omega_{z,\eps}}&&\leq c  \eps^{d/2} \left\| \nabla \psi^0 (z)\cdot{\bf n }\right\|_{-1/2,\partial \omega},\\\label{est-2}
\left\| \varphi^\eps\right\|_{1,\Omega_R} &&\leq c  \eps^{d/2} \left\| \nabla \psi^0 (z)\cdot{\bf n }\right\|_{-1/2,\partial \omega}.
\end{eqnarray}
The last estimates follow from a change of variable and the integral representation of the function $\varphi$. For more details and simular proof, one can consult \cite{GGM} for the elasticity problem or \cite{BHJM} for the Stokes problem. 

Now, we examine the boundary condition satisfied by $\phi^\eps$ on $\partial \O_{z,\eps}$. Using the smoothness of $\gamma$ and $\nabla\psi^0$ near the point $z$, one can derive
\begin{eqnarray}\label{est-3}
\left\| \gamma(x) [\nabla \psi^0(x)- \nabla \psi^0(z)]\cdot{\bf n }\right\|_{-1/2,\partial \O_{z,\eps}} \leq c\, \eps^{d/2}.
\end{eqnarray}
Finally, combining (\ref{est-1}), (\ref{est-2}) and (\ref{est-3}), one can deduce that the function $\phi^\eps=\psi^\eps-\psi^0-\eps \varphi((x-z)/\eps)$ satisfied the desired estimate
\begin{eqnarray*}
\left\| \phi^\eps \right\|_{1,\Omega_{z,\eps}}\leq c\,\eps^{d/2}.
\end{eqnarray*}

\section{Variation of the function $\K$}\label{cost-funct-variation}
This section is focused on the variation of the Kohn-Vogelius function $\K$ with respect to the presence of the small object $\O_{z,\eps}$ inside the domain $\Omega$. We will derive a simplified expression of the variation $\K (\Omega_{z,\eps})-\K (\Omega)$. The obtained result is presented in the following theorem.

\begin{theorem} \label{K-variation}  Let $\O_{z,\eps}$ be an arbitrary shaped object inserted inside the background domain $\Omega$. The variation $\K (\Omega_{z,\eps})-\K (\Omega)$ admits the expression
\begin{eqnarray}\label{decomp-k}
\nonumber
\K (\Omega_{z,\eps})-\K (\Omega) &=& \int_{\partial\O_{z,\eps}}\gamma(x) \nabla \psi^0_D \cdot{\bf n } \, \psi^\eps_D\,ds  -\int_{\partial\O_{z,\eps}}\gamma(x) \nabla \psi^0_N  \cdot{\bf n } \, \psi^\eps_N\,ds \\
&& + \int_{\O_{z,\eps}} F\,[\psi^0_D -\psi^0_N ]\,dx.
\end{eqnarray}
\end{theorem}

\noindent {\it Proof:} From the definition of $\K$, we have
\begin{eqnarray*}
\K (\Omega_{z,\eps})-\K (\Omega) &=& \int_{\Omega_{z,\eps}}\gamma(x)|\nabla \psi^\eps_D - \nabla \psi^\eps_N|^2\,dx- \int_{\Omega}\gamma(x)|\nabla \psi^0_D - \nabla \psi^0_N|^2\,dx\\
&=& \int_{\Omega_{z,\eps}}\gamma(x)| \nabla \psi^\eps_N|^2\,dx + \int_{\Omega_{z,\eps}}\gamma(x)|\nabla \psi^\eps_D|^2\,dx -2 \int_{\Omega_{z,\eps}}\gamma(x)\nabla \psi^\eps_D . \nabla \psi^\eps_N\,dx\\
&& -\int_{\Omega}\gamma(x)| \nabla \psi^0_N|^2\,dx - \int_{\Omega}\gamma(x)|\nabla \psi^0_D|^2\,dx -2 \int_{\Omega}\gamma(x)\nabla \psi^0_D . \nabla \psi^0_N\,dx.
\end{eqnarray*}
Then, the variation $\K (\Omega_{z,\eps})-\K (\Omega)$ can be decomposed as
$$
\K (\Omega_{z,\eps})-\K (\Omega)= T_N(\eps)+T_D(\eps)-2 T_M(\eps),
$$
where $T_N$ is the Neumann term 
$$
T_N(\eps)= \int_{\Omega_{z,\eps}}\gamma(x)| \nabla \psi^\eps_N|^2\,dx - \int_{\Omega}\gamma(x)| \nabla \psi^0_N|^2\,dx,
$$
$T_D$ is the Dirichlet  term 
$$
T_D(\eps) = \int_{\Omega_{z,\eps}}\gamma(x)|\nabla \psi^\eps_D|^2\,dx - \int_{\Omega}\gamma(x)|\nabla \psi^0_D|^2\,dx,
$$
and $T_M$ is the mixed term
$$
T_M (\eps)= \int_{\Omega_{z,\eps}}\gamma(x)\nabla \psi^\eps_D . \nabla \psi^\eps_N\,dx - \int_{\Omega}\gamma(x)\nabla \psi^0_D . \nabla \psi^0_N\,dx.
$$
Next, we shall examine each term separately. \\

\noindent $-$ {\it Calculate of the Neumann term.}

 \noindent From the weak formulation of the problems $ (\PB^\eps_N)$ and $ (\PB^0_N)$ one can obtain
\begin{eqnarray*}
T_N(\eps)&=& \int_{\Omega_{z,\eps}}F\,\psi^\eps_N\,dx + \int_{\Gamma_a}\Phi\,\psi^\eps_N \hbox{ds}  - \int_{\Omega}F\,\psi^0_N\,dx - \int_{\Gamma_a}\Phi\,\psi^0_N \hbox{ds}\\
&=& \int_{\Omega_{z,\eps}}F (\psi^\eps_N -  \psi^0_N)\,dx + \int_{\Gamma_a}\Phi (\psi^\eps_N -  \psi^0_N) \hbox{ds} - \int_{\O_{z,\eps}}F\,\psi^0_N\,dx.
\end{eqnarray*}
Then, it follows
\begin{eqnarray*}
T_N(\eps)&=& \int_{\Omega_{z,\eps}}\nabla \psi^\eps_N .\nabla (\psi^\eps_N -  \psi^0_N)\,dx  - \int_{\O_{z,\eps}}F\,\psi^0_N\,dx.
\end{eqnarray*}
Using the Green formula and the fact that 
$$\mbox{ div } (\gamma(x) \nabla (\psi^\eps_N -  \psi^0_N))  =0 \hbox{ in } \Omega_{z,\eps},\,\, \psi^\eps_N -  \psi^0_N=0 \hbox{ on }\Gamma_i\, \hbox{ and }\,  \gamma(x) \nabla (\psi^\eps_N -  \psi^0_N)  \cdot{\bf n }  = 0 \hbox{ on }\Gamma_a,
$$ 
we deduce
\begin{eqnarray}\nonumber
T_N(\eps)&=& \int_{\partial\O_{z,\eps}} \gamma(x) \nabla (\psi^\eps_N -  \psi^0_N) \cdot{\bf n }  \psi^\eps_N \,ds  - \int_{\O_{z,\eps}}F\,\psi^0_N\,dx,\\
\label{term-neumann}
&=& - \int_{\partial\O_{z,\eps}} \gamma(x) \nabla  \psi^0_N  \cdot{\bf n }  \psi^\eps_N \,ds  - \int_{\O_{z,\eps}}F\,\psi^0_N\,dx.
\end{eqnarray}

\noindent $-$ {\it Calculate of the Dirichlet term.}

\noindent  We have 
\begin{eqnarray*}
T_D(\eps) &=& \int_{\Omega_{z,\eps}}\gamma(x)|\nabla \psi^\eps_D|^2\,dx - \int_{\Omega}\gamma(x)|\nabla \psi^0_D|^2\,dx,\\
&=& \int_{\Omega_{z,\eps}}\gamma(x) \nabla (\psi^\eps_D -  \psi^0_D)\nabla (\psi^\eps_D +  \psi^0_D)\,dx - \int_{\O_{z,\eps}}\gamma(x)|\nabla \psi^0_D|^2\,dx.
\end{eqnarray*}
Using the Green formula and the fact that 
$$\mbox{ div } (\gamma(x) \nabla (\psi^\eps_D +  \psi^0_D))  =2F \hbox{ in } \Omega_{z,\eps} \quad \hbox{ and }\quad  \psi^\eps_D -  \psi^0_D  = 0 \hbox{ on }\Gamma,
$$ 
we derive
\begin{eqnarray*}
T_D(\eps)&=& \int_{\partial\O_{z,\eps}} \gamma(x) \nabla (\psi^\eps_D + \psi^0_D) \cdot{\bf n }   (\psi^\eps_D - \psi^0_D)\,ds + 2 \int_{\Omega_{z,\eps}}F (\psi^\eps_D -  \psi^0_D)\,dx \\
&&- \int_{\O_{z,\eps}}\gamma(x)|\nabla \psi^0_D|^2\,dx.
\end{eqnarray*}
Recalling that $\mbox{ div } (\gamma(x) \nabla \psi^0_D)  =F \hbox{ in } \O_{z,\eps}$ and taking into account of the normal orientation one can write
$$
\int_{\O_{z,\eps}}\gamma(x)|\nabla \psi^0_D|^2\,dx= \int_{\O_{z,\eps}} F \, \psi^0_D\,dx - \int_{\partial\O_{z,\eps}} \gamma(x) \nabla \psi^0_D  \cdot{\bf n } \, \psi^0_D \,ds.
$$
Then, we obtain
\begin{eqnarray}\label{term-dirichlet}
T_D(\eps)&=& \int_{\partial\O_{z,\eps}} \gamma(x) \nabla \psi^0_D  \cdot{\bf n }   \psi^\eps_D \,ds + 2 \int_{\Omega_{z,\eps}}F (\psi^\eps_D -  \psi^0_D)\,dx - \int_{\O_{z,\eps}} F \, \psi^0_D\,dx.
\end{eqnarray}

\noindent $-$ {\it Calculate of the mixed term.}

\noindent From the weak formulation of the problems $ (\PB^\eps_N)$ and $ (\PB^0_N)$ one can derive
\begin{eqnarray*}
T_M (\eps)&=& \int_{\Omega_{z,\eps}}\gamma(x)\nabla \psi^\eps_D . \nabla \psi^\eps_N\,dx - \int_{\Omega}\gamma(x)\nabla \psi^0_D . \nabla \psi^0_N\,dx\\
&=&\int_{\Omega_{z,\eps}}F\,\psi^\eps_D\,dx + \int_{\Gamma_a}\Phi\,\psi^\eps_D \hbox{ds} -\int_{\Omega}F\,\psi^0_D\,dx - \int_{\Gamma_a}\Phi\,\psi^0_D \hbox{ds} 
\end{eqnarray*}
Using the fact that $\psi^\eps_D -  \psi^0_D  = 0 \hbox{ on }\Gamma$, it follows
\begin{eqnarray}\label{term-mixed}
T_M (\eps)&=& \int_{\Omega_{z,\eps}}F \,(\psi^\eps_D -  \psi^0_D)\,dx -\int_{\O_{z,\eps}}F\,\psi^0_D\,dx.
\end{eqnarray}

\noindent Exploiting the obtained  expressions (\ref{term-neumann}), (\ref{term-dirichlet}) and  (\ref{term-mixed}), one can easily deduce the desired result of Theorem \ref{K-variation}.


\section{Asymptotic expansion}\label{asymptotic}
In this section, we derive a topological asymptotic expansion for the Kohn-Vogelius function $\K$. The mathematical analysis is general and can be adapted for various partial differential equations. 

To this end, we introduce the polarization matrix $\mathcal M_\omega$. Thanks to the linearity of the integral equations (\ref{integ-N}) and (\ref{integ-D}), there exists a $d\times d$ matrix $\mathcal M_\omega$ such that
$$
\dint_{\partial \omega} \eta_N (y) y^T ds(y)=\mathcal M_\omega \nabla \psi^0_N(z) \mbox{ and } \dint_{\partial \omega} \eta_D (y) y^T ds(y)=\mathcal M_\omega \nabla \psi^0_D(z).
$$
The matrix $\mathcal M_\omega$ can be defined as
$$
(\mathcal M_\omega)_{ij}=\int_{\partial \omega} \eta_i y_j ds(y),\quad 1\leq i,j\leq d,
$$
where $y_j$ is the $j$th component of $y\in \mathbb{R}^d$ and $\eta_i $ is the solution to 
$$
-\frac{\eta_i(y)}{2}+\int_{\partial\omega}\nabla E(y-x)   \cdot{\bf n } \,\eta_i(x)ds(x)=-e_i   \cdot{\bf n } ,\quad y\in \partial \omega
$$
with $\{e_i)_{1\leq i\leq d}$ is the canonical basis in $\mathbb{R}^d$.

The topological sensitivity analysis with respect to the presence of an arbitrary shaped object is described by the following Theorem. 

\begin{theorem}\label{th-kv}  Let $\O_{z,\eps}$ be an arbitrary shaped object inserted inside the background domain $\Omega$. The function $\K$ admits the asymptotic expansion
$$
\K(\Omega\backslash\overline{\O_{z,\eps}})=\K (\Omega)+ \eps^d \delta \K(z) + o(\eps^d),
$$
with $\delta \K$ is the topological gradient 
\begin{eqnarray*}
\delta \K(x) &=& \gamma(x) \left\{\nabla \psi^0_N(x) \mathcal M_\omega \nabla \psi^0_N(z)- \nabla \psi^0_D(x) \mathcal M_\omega \nabla \psi^0_D(z) \right\}\\
&& -2 |\omega|\, F(x) (\psi^0_N(x)-\psi^0_D(x)),\,\forall x\in \Omega.
\end{eqnarray*}
\end{theorem}
The polarization matrix $\mathcal M_\omega$ can be determined analytically in some cases. Otherwise, it can be approximated numerically. \\
Particularly, in the case of circular or spherical object (i.e. $\omega = B(0,1)$), the matrix $\mathcal M_\omega$ is given by
$$
\mathcal M_\omega = 2\pi \mathcal I \mbox{ if } d=2 \mbox{ or } d=3,
$$
where $\mathcal I $ is the $d\times d$ identity matrix.

\begin{corollary}\label{cor-bal} If $\omega = B(0,1)$, the function $\K$ satisfies the following asymptotic expansion
\begin{eqnarray*}
 \K(\Omega_{z,\eps})-\K (\Omega)=2 \pi \,\varepsilon^d \delta \K(z)+o(\varepsilon^d),
\end{eqnarray*}
and the topological gradient $\delta \K$ admits the expression
\begin{eqnarray*}
\delta \K(x) = \left\{
\begin{array}{ll}
\gamma(x)\left(\left|\nabla \psi^0_N(x)\right|^2-\left|\nabla \psi^0_D(x)\right|^2\right)-F(x)\left( \psi^0_N(x)-\psi^0_D(x)\right),&\mbox{ if } d=2,\\
\gamma(x)\left(\left|\nabla \psi^0_N(x)\right|^2-\left|\nabla \psi^0_D(x)\right|^2\right)-\displaystyle \frac{4}{3}F(x)\left( \psi^0_N(x)-\psi^0_D(x)\right),&\mbox{ if } d=3.
\end{array}
\right.
\end{eqnarray*}
\end{corollary} 

\noindent {\it Proof of Theorem \ref{th-kv}} : 
It is established in Theorem \ref{K-variation} that  the variation of the function $\K$  can be rewritten as
\begin{eqnarray}\label{decomp-k}
\nonumber
\K (\Omega_{z,\eps})-\K (\Omega) &=& \int_{\partial\O_{z,\eps}}\gamma(x) \nabla \psi^0_D\cdot{\bf n }  \, \psi^\eps_D\,ds  -\int_{\partial\O_{z,\eps}}\gamma(x) \nabla \psi^0_N \cdot{\bf n }  \, \psi^\eps_N\,ds \\
&& + \int_{\O_{z,\eps}} F\,(\psi^0_D -\psi^0_N )\,dx.
\end{eqnarray}
Next, we will derive an asymptotic expansion with respect to $\eps$ for each term separately.

\vspace{0.25cm}
\noindent $-$ {\it Asymptotic expansion for the first integral term}.

\noindent The first term in (\ref{decomp-k}) can be decomposed as 
\begin{eqnarray*}
\int_{\partial\O_{z,\eps}}\gamma(x) \nabla \psi^0_D\cdot{\bf n }  \, \psi^\eps_D\,ds = \int_{\partial\O_{z,\eps}}\gamma(x) \nabla \psi^0_D\cdot{\bf n }  \, (\psi^\eps_D(x)-\psi^0_D(x)-\eps \varphi_D((x-z)/\eps))\,ds\\
+ \int_{\partial\O_{z,\eps}}\gamma(x) \nabla \psi^0_D\cdot{\bf n }  \, \psi^0_D\,ds+ \eps \int_{\partial\O_{z,\eps}}\gamma(x) \nabla \psi^0_D\cdot{\bf n }  \, \varphi_D((x-z)/\eps)\,ds.
\end{eqnarray*}
Using Theorem \ref{est-pert-sol} and the fact that $\gamma(x)\nabla \psi^0_D(x)$ is uniformly bounded in $\overline{\O}_{z,\eps}$, one can obtain
\begin{eqnarray}\label{estim-1}
\int_{\partial\O_{z,\eps}}\gamma(x) \nabla \psi^0_D\cdot{\bf n }  \, [\psi^\eps_D(x)-\psi^0_D(x)-\eps \varphi_D((x-z)/\eps)]\,ds  = o(\eps^d).
\end{eqnarray}
 Due to the Green's formula and taking into account of the normal orientation,
\begin{eqnarray}\label{estim-2}
\int_{\partial\O_{z,\eps}}\gamma(x) \nabla \psi^0_D\cdot{\bf n } (x) \, \psi^0_D(x)ds=-\int_{\O_{z,\eps}}\gamma(x) \nabla \psi^0_D(x) .\nabla \psi^0_D(x) dx + \int_{\O_{z,\eps}} F(x)\psi^0_D(x) dx.
\end{eqnarray}
By the change of variable $x=z+\eps y$, one can write
\begin{eqnarray*}
\eps \int_{\partial\O_{z,\eps}}\gamma(x) \nabla \psi^0_D\cdot{\bf n }  \, \varphi_D((x-z)/\eps)ds= \eps^d \int_{\partial\omega}\gamma(z) \nabla \psi^0_D(z)\cdot{\bf n }  \, \varphi_D(y)ds(y)\\
+ \eps^d \int_{\partial\omega}[\gamma(z+\eps y)\nabla \psi^0_D(z+\eps y)-\gamma(z) \nabla \psi^0_D(z)]\cdot{\bf n }  \, \varphi_D(y)ds(y).
\end{eqnarray*}
Using again the smoothness of the function $x\mapsto \gamma(x)\nabla \psi^0_D(x)$ near $z$, one can obtain
\begin{eqnarray*}
\eps \int_{\partial\O_{z,\eps}}\gamma(x) \nabla \psi^0_D\cdot{\bf n }  \, \varphi_D((x-z)/\eps)ds= \eps^d \gamma(z) \int_{\partial\omega} \nabla \psi^0_D(z)\cdot{\bf n }  \, \varphi_D(y)ds(y)+o(\eps^d).
\end{eqnarray*}
Let $\xi^{\varphi_D}$ the unique solution to
 \begin{equation*}
\left\{
\begin{array}{rll}
\displaystyle
- \Delta \xi^{\varphi_D}  =&  0 &\hbox{ in }  \omega \\
\xi^{\varphi_D} =& \varphi_D & \hbox{ on } \partial \omega .
\end{array}
\right.
\end{equation*}
Exploiting  the Green's formula, the fact that $\Delta (\nabla \psi^0_D(z)y)=0$ in $\omega$ and taking into account of the normal orientation, one can derive 
\begin{eqnarray*}
\int_{\partial\omega} \nabla \psi^0_D(z)\cdot{\bf n }  \, \varphi_D(y)ds(y)&=& - \int_{\omega} \nabla \psi^0_D(z). \nabla \xi^{\varphi_D} (y)dy\\
&=& \int_{\partial\omega} \nabla \xi^{\varphi_D} (y)\cdot{\bf n }  [\nabla \psi^0_D(z)y]ds(y).
\end{eqnarray*}
Recalling that the density $\eta_D$ is defined as the jump of the flux through the boundary $\partial \omega$
\begin{eqnarray*}
\eta_D(y)= -\nabla \psi^0_D(z)\cdot{\bf n }  - \nabla \xi^{\varphi_D} (y)\cdot{\bf n } ,\quad y\in \partial \omega.
\end{eqnarray*}
It follows 
\begin{eqnarray*}
\int_{\partial\omega} \nabla \psi^0_D(z)\cdot{\bf n }  \, \varphi_D(y)ds(y)=- \int_{\omega} [\nabla \psi^0_D(z)y]\eta_D(y)dy - \int_{\partial\omega}\nabla \psi^0_D(z)\cdot{\bf n } [\nabla \psi^0_D(z)y]ds(y).
\end{eqnarray*}
Then, we deduce
\begin{eqnarray}\label{estim-3}
\nonumber
\eps \int_{\partial\O_{z,\eps}}\gamma(x) \nabla \psi^0_D\cdot{\bf n }  \, \varphi_D((x-z)/\eps)ds= -\eps^d \gamma(z) \nabla \psi^0_D(x) .\dint_{\partial \omega} \eta_D (y) y ds(y)\\
+\int_{\O_{z,\eps}}\gamma(z)\nabla \psi^0_D(z).\nabla \psi^0_D(z)dx+o(\eps^d).
\end{eqnarray}
Combining (\ref{estim-1}), (\ref{estim-2}) and (\ref{estim-3}), we have  
\begin{eqnarray*}
\int_{\partial\O_{z,\eps}}\gamma(x) \nabla \psi^0_D\cdot{\bf n }  \, \psi^\eps_D\,ds =  -\eps^d \gamma(z) \nabla \psi^0_D(x) .\dint_{\partial \omega} \eta_D (y) y ds(y) + \int_{\O_{z,\eps}} F(x)\psi^0_D(x) dx\\
+\int_{\O_{z,\eps}}\gamma(z)\nabla \psi^0_D(z).\nabla \psi^0_D(z)dy-\int_{\O_{z,\eps}}\gamma(x) \nabla \psi^0_D(x) .\nabla \psi^0_D(x) dx +o(\eps^d).
\end{eqnarray*}
Finally, using a change of variable and the Taylor's formula one can easily deduce that the first term in (\ref{decomp-k}) admits the following expansion
\begin{eqnarray}\label{est-term1}
\int_{\partial\O_{z,\eps}}\gamma(x) \nabla \psi^0_D\cdot{\bf n }  \, \psi^\eps_D\,ds =-\eps^d \gamma(z) \nabla \psi^0_D(x) .\dint_{\partial \omega} \eta_D (y) y ds(y) + \eps^d |\omega|\, F(z)\psi^0_D(z)+o(\eps^d).
\end{eqnarray}

\vspace{0.25cm}
\noindent $-$ {\it Asymptotic expansion for the second integral term}.

\noindent \noindent The second term in (\ref{decomp-k}) can be decomposed as 
\begin{eqnarray*}
\int_{\partial\O_{z,\eps}}\gamma(x) \nabla \psi^0_N\cdot{\bf n }  \, \psi^\eps_N\,ds = \int_{\partial\O_{z,\eps}}\gamma(x) \nabla \psi^0_N\cdot{\bf n }  (x) (\psi^\eps_N(x)-\psi^0_N(x)-\eps \varphi_N((x-z)/\eps))\,ds\\
+ \int_{\partial\O_{z,\eps}}\gamma(x) \nabla \psi^0_N(x)\cdot{\bf n }  \, \psi^0_N(x)\,ds+ \eps \int_{\partial\O_{z,\eps}}\gamma(x) \nabla \psi^0_N\cdot{\bf n }  \, \varphi_N((x-z)/\eps)\,ds.
\end{eqnarray*}
The established estimates for the first term can be easily adapted for this term using $\psi^0_N$ instead of $\psi^0_D$ and $\varphi_N$ instead of $\varphi_D$. Then, one can prove that the second integral term in (\ref{decomp-k}) satisfies the following expansion
\begin{eqnarray}\label{est-term2}
\int_{\partial\O_{z,\eps}}\gamma(x) \nabla \psi^0_N\cdot{\bf n }  \, \psi^\eps_N \,ds =-\eps^d \gamma(z) \nabla \psi^0_N(x) .\dint_{\partial \omega} \eta_N (y) y ds(y) + \eps^d |\omega|\, F(z)\psi^0_N(z)+o(\eps^d).
\end{eqnarray}

\vspace{0.25cm}
\noindent $-$ {\it Asymptotic expansion for the third integral term}.

\noindent The estimate of the third term is based on the change of variable $y=z+ \eps y$ and the smoothness of $\gamma$ in $\O_{z,\eps}$. Using Taylor's Theorem, one can derive
\begin{eqnarray}\label{est-term3}
\int_{\O_{z,\eps}} F\,(\psi^0_D -\psi^0_N )dx = \eps^d |\omega| F(z) (\psi^0_D(z) -\psi^0_N (z)) + o(\eps^d).
\end{eqnarray}
The desired asymptotic expansion of the function $\K$ follows from the estimates  (\ref{est-term1}),  (\ref{est-term2}) and  (\ref{est-term3}).

\section{Algorithm and numerical results}
In this section we consider the bidimentional case and we present a fast and simple one-iteration identification algorithm. Our numerical procedure is based on the formula described by  Corollary \ref{cor-bal}.

The unknow object $\O$ is identified using a level set curve of the topological gradient $\delta \K$. More precisely, the unknow object $\O$ is likely to be located at zone where the topological gradient $\delta \K$ is negative.

\vskip 0.4cm

\noindent{\bf One-iteration algorithm:}
\begin{itemize}
\item Solve the two problems ($\PB_N^0$) and ($\PB_D^0$),
\item Compute the topological gradient $\delta \K(x),\, x\in \Omega$,
\item Determine the unknow object 
$$\O =\left\{x\in \Omega;\hbox{ such that } \delta \K(x)<c<0\right\},
$$
where $c$ is a constant chosen in such a way that the cost function $\K$ decreases as most as possible.
\end{itemize}
\vskip 0.2cm
This one-iteration procedure has already been illustrated in \cite{AHM2} for the identification of cracks from overdetermined boundary data and  in \cite{BHJM} for the detection of small gas bubbles in Stokes flow.

 Next, we presente some numerical results showing the efficiency and accuracy of our proposed one-iteration algorithm. In Figure 1, we test our algorithm on circular shape. In Figure 2, we consider the case of an elliptical shape. As one can observe, the domain to be detected is located at zone where the topological gradient is negative and it is approximated by a level set curve of the topological gradient $\delta \K$. The result is quite efficient. In Figure 3, we obtain an intersting reconstruction result for a non trivial shape.


\begin{figure}[h]
\begin{center}
\includegraphics[width=4.5cm, height=3.1cm]{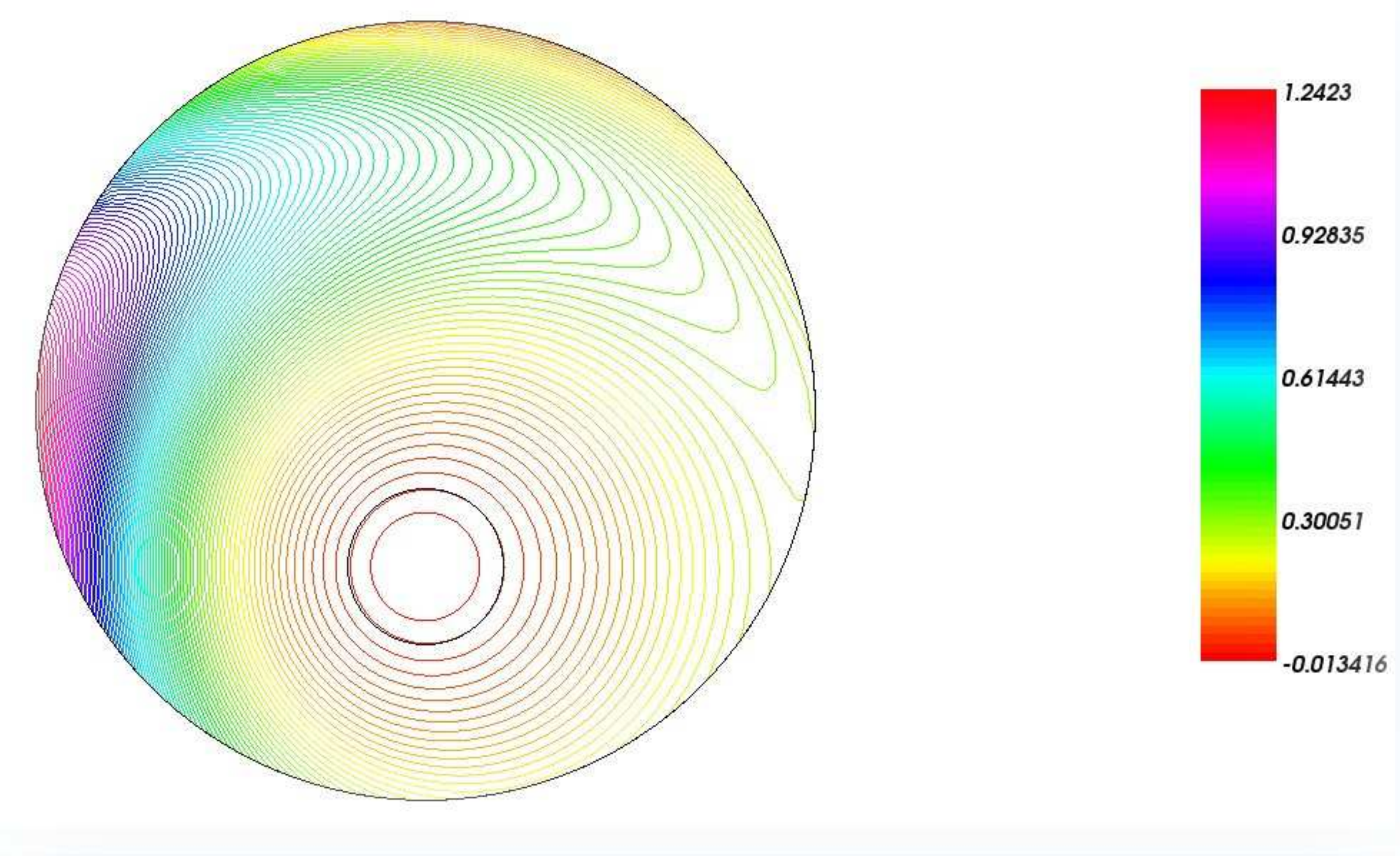}
\includegraphics[width=4.5cm, height=3.1cm]{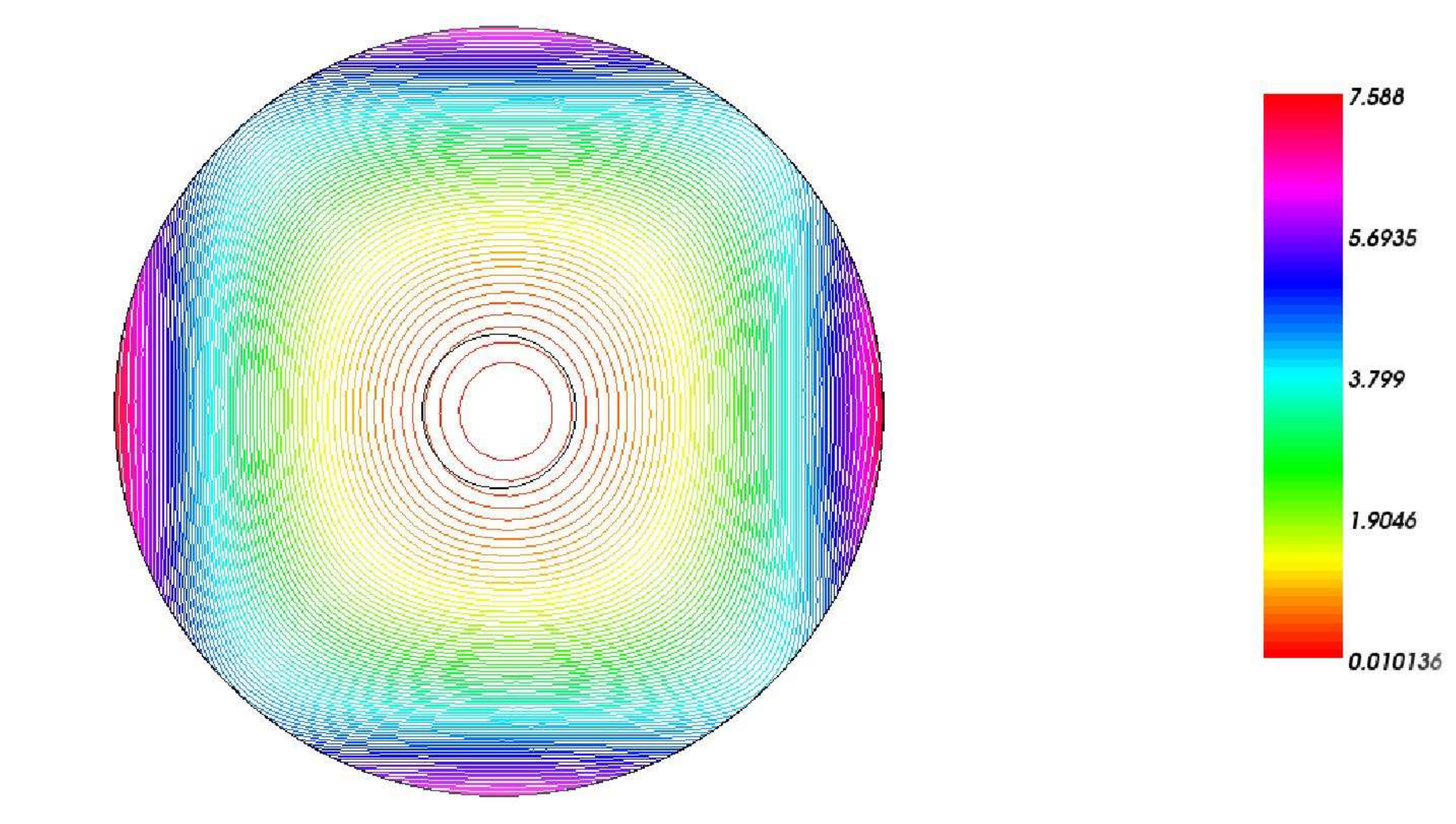}
\includegraphics[width=4.5cm, height=3.1cm]{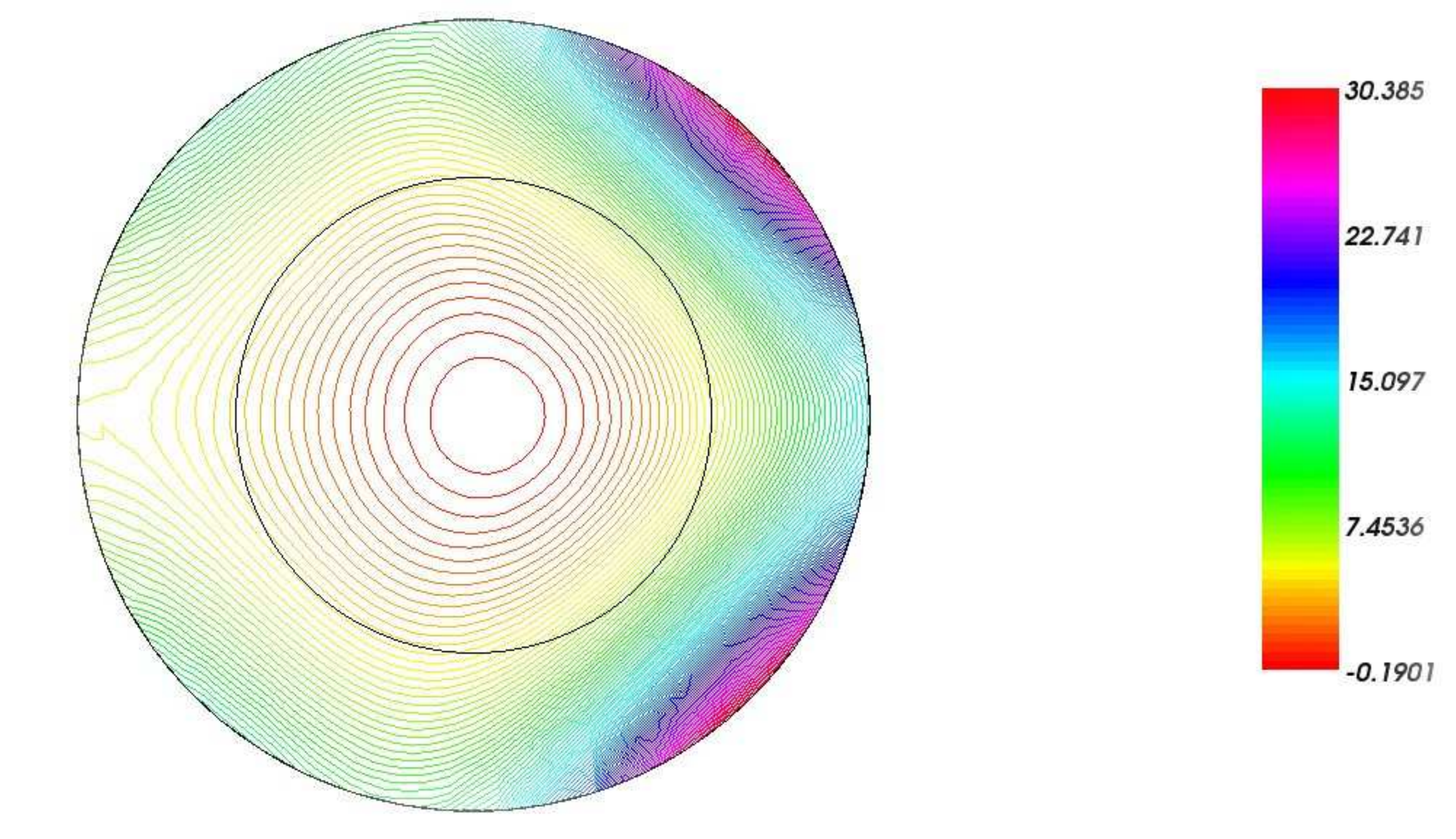}
\end{center}

\vskip -0.6cm
\caption{Reconstruction of circle shaped objects}
\end{figure}

\begin{figure}[h]
\begin{center}
\includegraphics[width=4.5cm, height= 3.1cm]{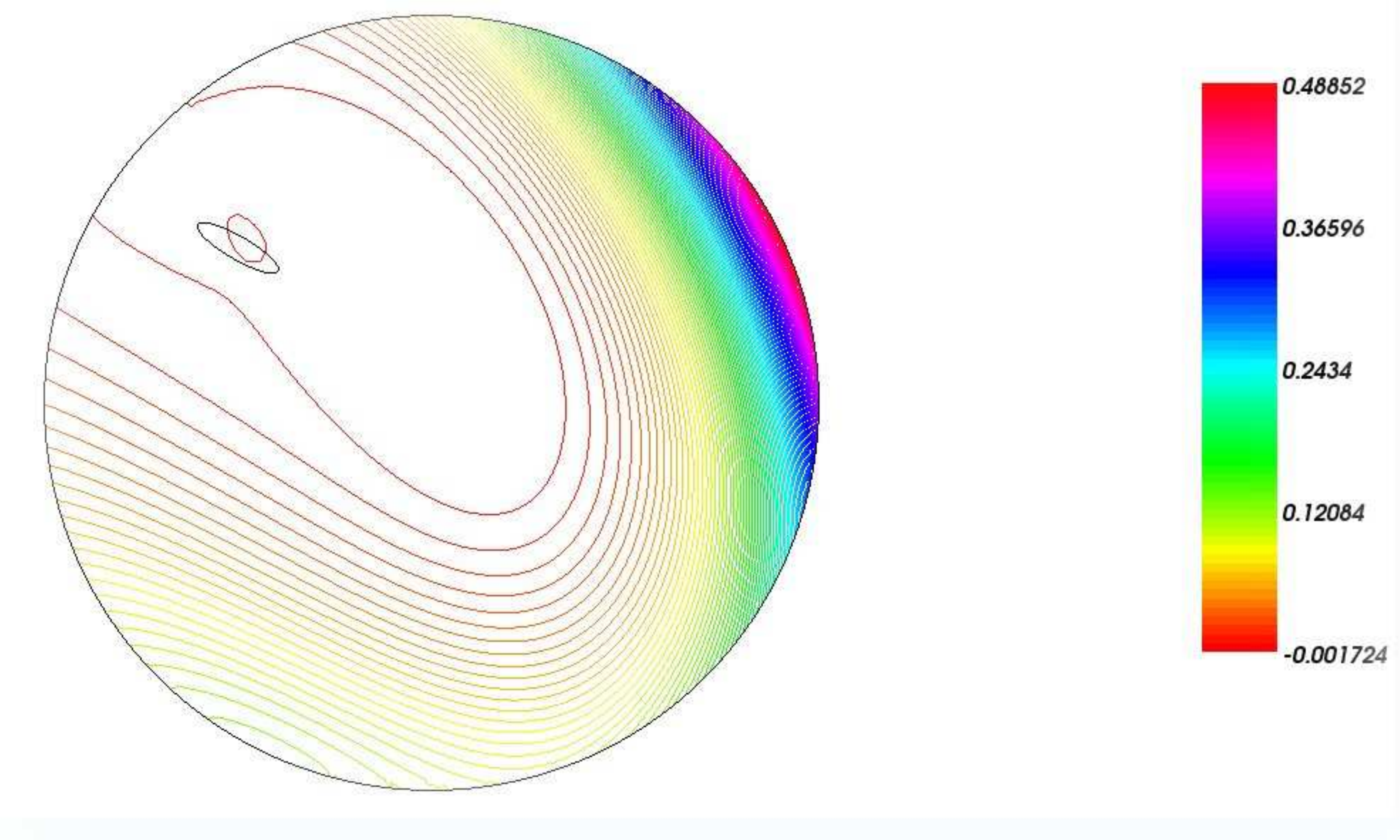}
\includegraphics[width=4.5cm, height= 3.cm]{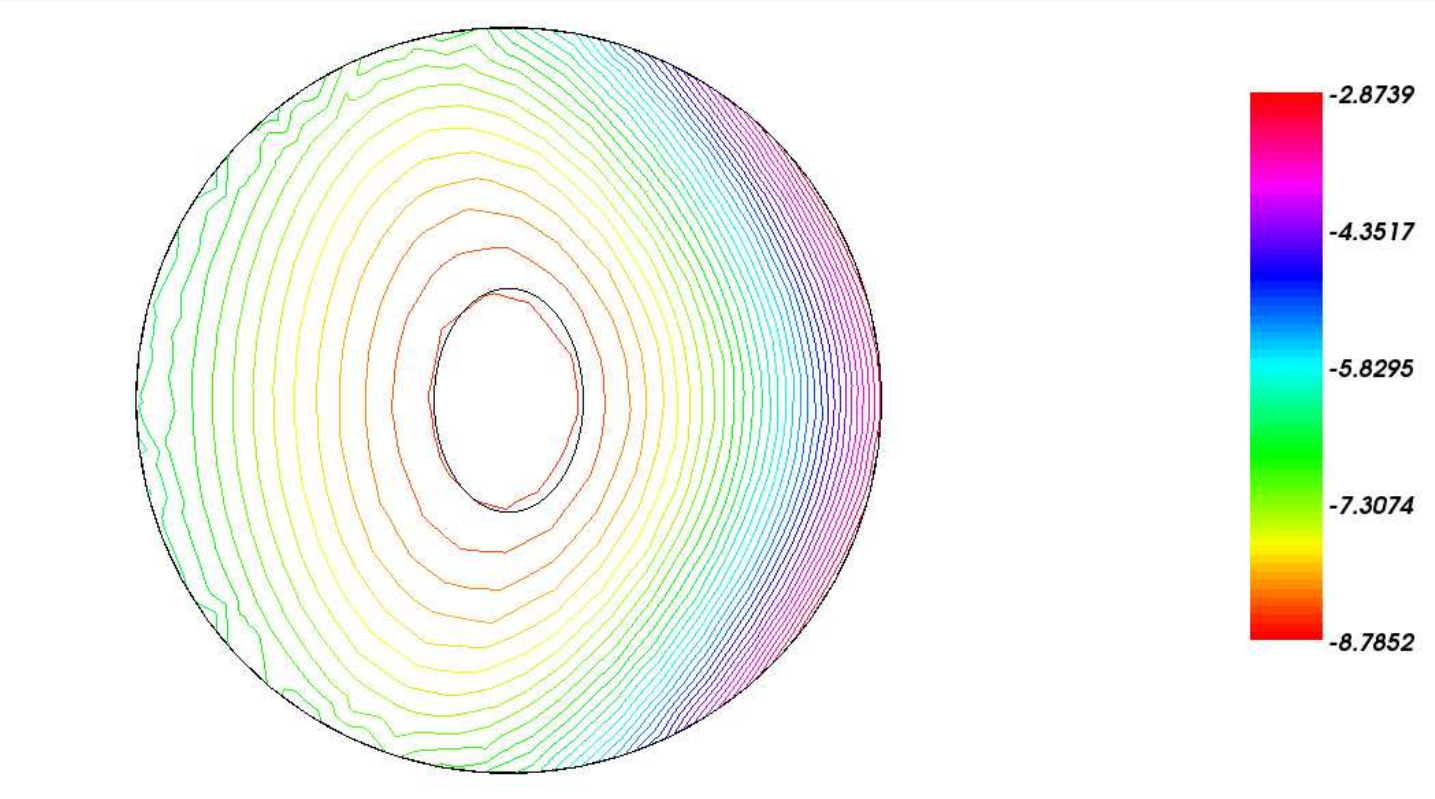}
\includegraphics[width=4.5cm, height= 3.1cm]{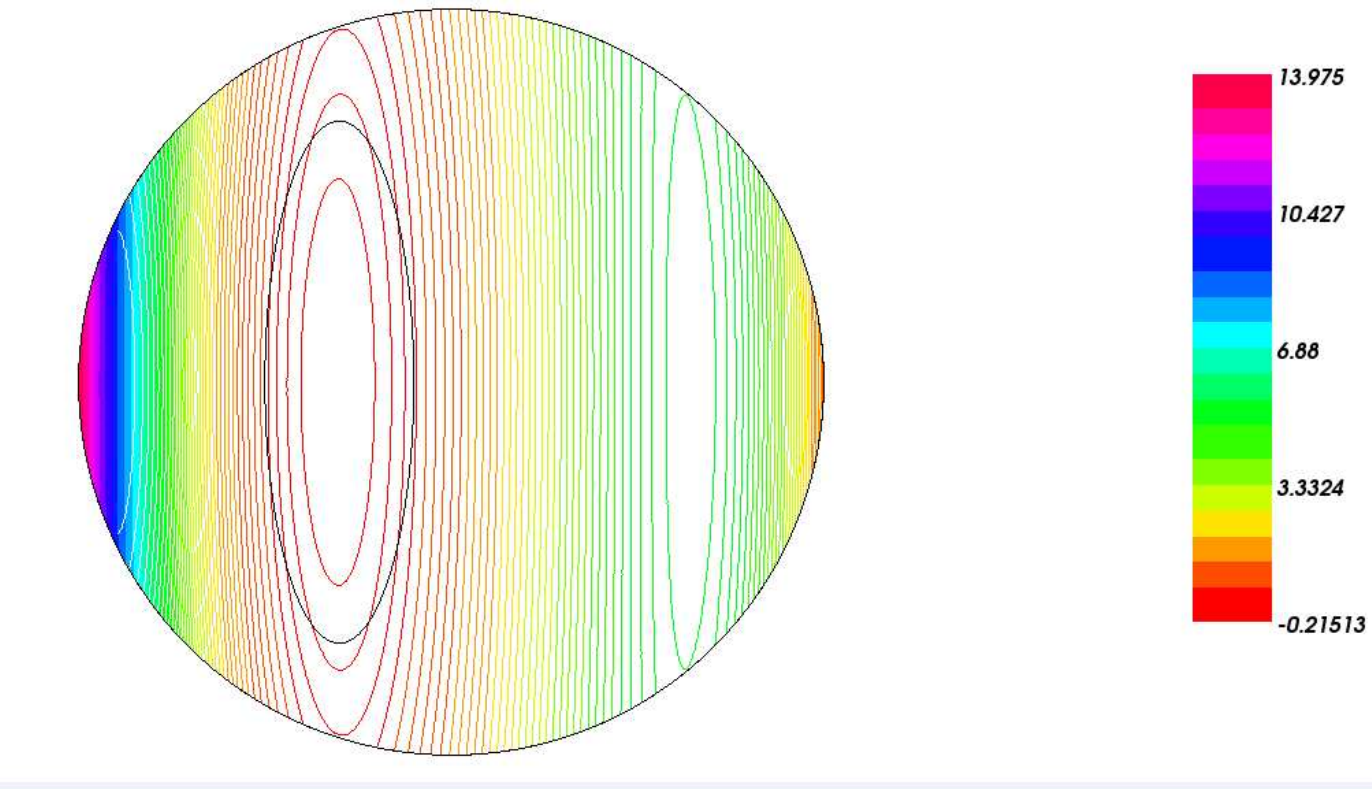}
\end{center}

\vskip -0.6cm
\caption{Reconstruction of an ellipse shaped objects}
\begin{center}
\includegraphics[width=5.cm, height= 3.1cm]{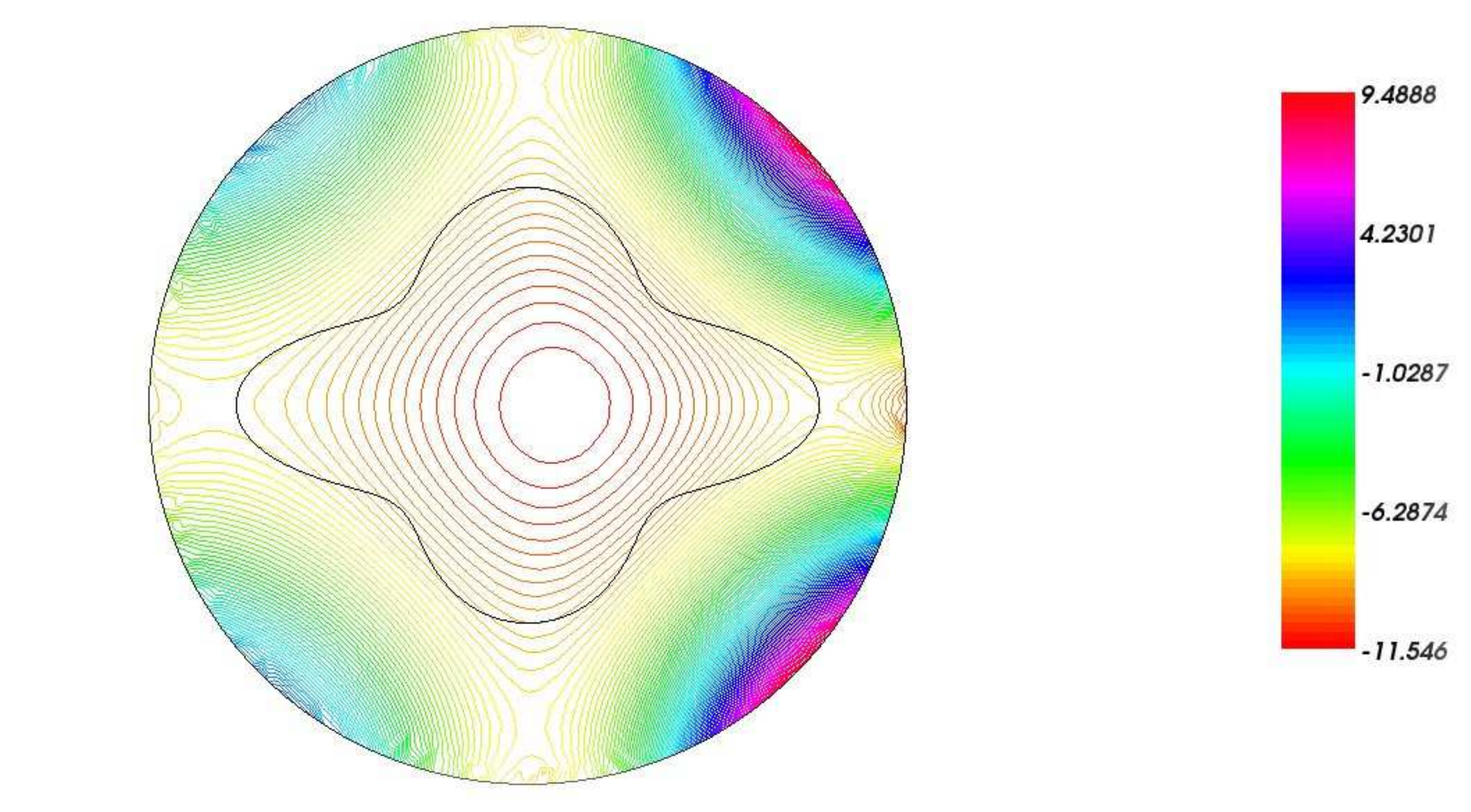}
\end{center}

\vskip -0.6cm
\caption{Reconstruction of a non trivial shape}
\end{figure}

In Figure 4 we illustrate the case of a complex geometry. As one can see, we detect efficiency the location of the unknown domain but not its shape. The obtained result can serve as a good initial guess for an iterative optimization process based on the shape derivative.
\begin{figure}[h]
\begin{center}
\includegraphics[width=5.cm, height= 3.cm]{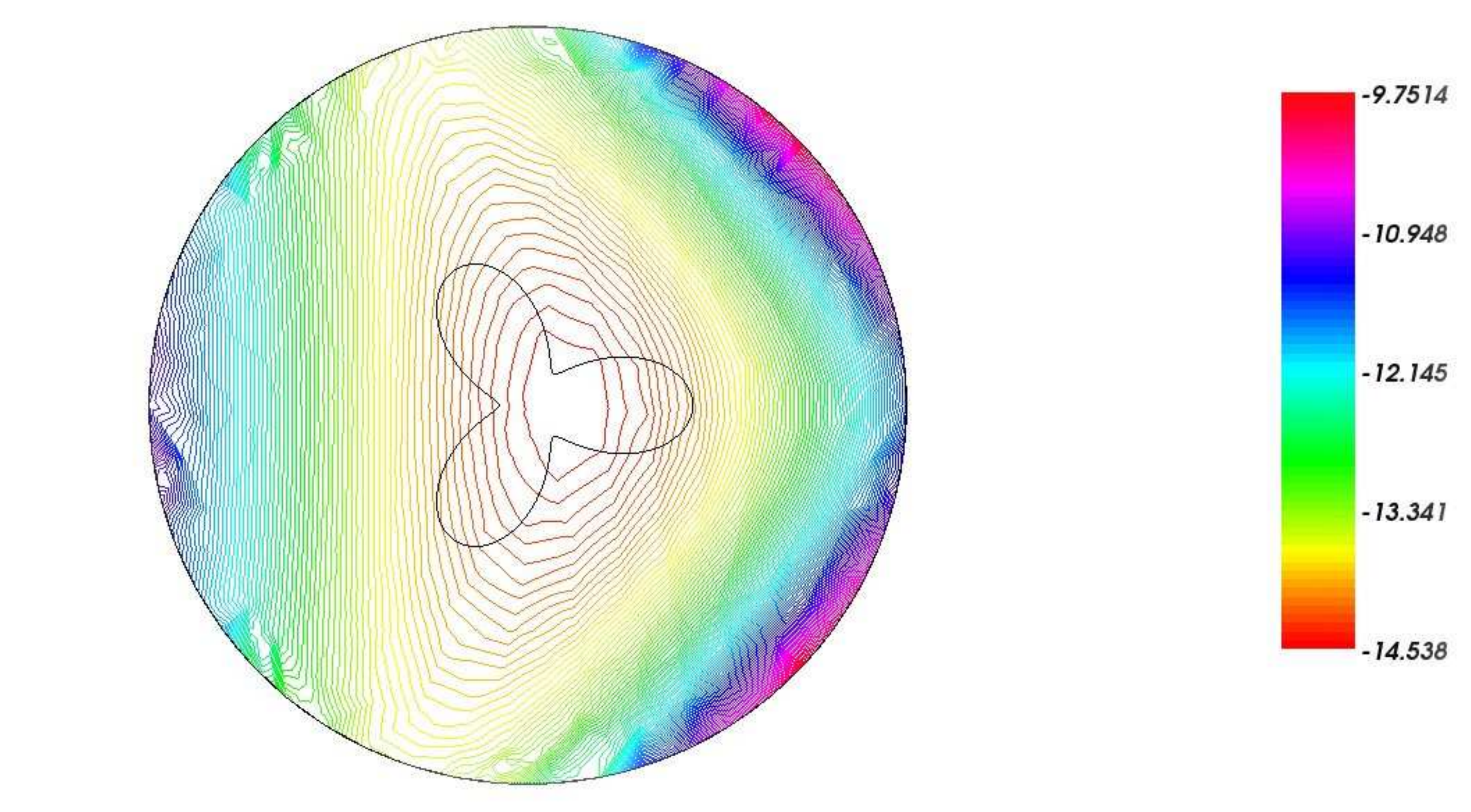}
\end{center}

\vskip -0.6cm
\caption{Reconstruction of a complex shape}

\vskip 0.3cm
\end{figure}

\section{Conclusion}
The presented work concerns the detection of objects immersed in anisotropic media from boundary measurements. The present approach  is based on the Kohn-Vogelius formulation and the topological gradient method.  A topological sensitivity analysis is derived for an energy like functional. An accurate and fast reconstruction algorithm is proposed. The efficiency and accuracy of the suggested  algorithm are illustrated by some numerical results.  The considered model can be viewed as a prototype of a geometric inverse problem arising in many applications. The presented approach is general and can be adapted for various partial differential equations.


\end{document}